\documentclass[12pt]{article}
\usepackage{txfonts}
\usepackage{amssymb}
\usepackage{amsfonts}
\usepackage{mathrsfs}
\newcounter{S}
 \makeatother

\newtheorem{definition}{{\bf Definition}}
\newtheorem{theorem}{{\bf Theorem}}
\newtheorem{corollary}{{\bf Corollary}}
\newtheorem{lemma}{{\bf Lemma}}
\newtheorem{example}{{\bf Example}}
\newtheorem{proposition}{\bf Proposition}
\newtheorem{remark}{\bf Remark}

\newcounter{A}
 \makeatother

 \makeatother

 \makeatother
 \makeatother
 \makeatother
 \makeatother

 \makeatother

\newenvironment{Definition}{\begin{definition}\hspace{-6pt}{\bf.}}{\end{definition}}
\newenvironment{Theorem}{\begin{theorem}\hspace{-6pt}{\bf.}}{\end{theorem}}
\newenvironment{Lemma}{ \begin{lemma}\hspace{-6pt}{\bf.}}{\end{lemma}}
\newenvironment{Proposition}{ \begin{proposition}\hspace{-6pt}{\bf.}}{\end{proposition}}
\newenvironment{Corollary}{\begin{corollary}\hspace{-6pt}{\bf.}}{\end{corollary}}
\newenvironment{Remark}{\begin{remark}\hspace{-6pt}{\bf.}}{\end{remark}}
\newenvironment{Example}{\begin{example}\hspace{-6pt}{\bf.}}{\end{example}}

\newcommand{\se}[1]{\addtocounter{S}{1}\setcounter{equation}{0}\setcounter{lemma}{0}\setcounter{proposition}{0}\setcounter{theorem}{0}\setcounter{remark}{0}\setcounter{corollary}{0}\setcounter{definition}{0}\section*{\bf\normalsize\arabic{S}.$\;$ #1}}

\begin{document}
\baselineskip 20pt
\parindent 0pt
\title{\Large\bf Second order subexponential distributions with finite mean and their applications to subordinated distributions}
\author{\bf Jianxi Lin\thanks{Author's Email address: linjx999@gmail.com}\\{\normalsize Mathematics School of Xiamen  University, Xiamen, Fujian 361005, China.}}
\maketitle
{\bf Abstract}\\
\mbox{}\quad
  Consider a probability distribution subordinate to a
  subexponential distribution with finite mean. In this paper, we discuss the second order tail behavior
  of the subordinated distribution within a rather general framework in which we do not require the existence of density functions.
 For this aim, the so-called second order subexponential distribution is proposed and some related properties of its are established.  Our results unified and improved some classical results.
 \se{Introduction}\label{dyjzyxz}\quad\quad
Let $N$ be a non-negative integer valued random variable with
distribution $\{ p_n \} _{n \ge 0} $ and $X_1$, $X_2$, $\cdots$ be a
sequence of non-negative i.i.d. random variables, independent of
$N$. The common distribution of $X_i$'s is denoted by $F$. Define for $n\ge 1$,
\begin{eqnarray}
S_n:=\sum\limits_{k=1}^{n}X_k,
\end{eqnarray}
and $S_0=0$. In many fields
of applied probability, one has to investigate the tail behavior of
$S_N$, whose distribution is equal to
\begin{eqnarray}\label{soGxs}
 G (x) := \sum\limits_{n = 0}^\infty  {p_n
{F^{n*} } (x)},
\end{eqnarray}
where $F^{*0}$ is the unit mass at zero and for $n\ge 1$, $F^{*n}$
denotes the $n$-fold convolution of distribution $F$. Obviously, $G$
is a probability distribution subordinate to $F$ with subordinator
$\{ p_n \} _{n \ge 0} $.
\\ \mbox{}\quad Denote the tail of distribution $G$  by
$\overline{G}(x)=G(x,\infty)=1-G(x)$. A first order approximation to
$\overline{G}(x)$ as $x\to \infty$ has been considered by Chistyakov
\cite{chistyakov1964}, in which he introduce the so-called
subexponential distribution class $\mathscr{S}$. By definition, a
distribution $F$ on $[0,\infty)$ is said to belong to the class
$\mathscr{S}$ if for $n=2$ (hence for all $n\ge 2$),
\begin{equation}\label{czsgx}
\mathop {\lim }\limits_{x \to \infty } \frac{{\overline {F^{n*} }
(x)}}{{\overline F (x)}} = n.
\end{equation}
 Chistyakov
\cite{chistyakov1964} states that if $F\in \mathscr{S}$ and
$E(z^{N})$ is analytic at $z=1$, then
\begin{equation}\label{sta}
\overline G (x)  \sim (\sum\limits_{n = 0}^\infty  {np_n }
)\overline F (x), \quad x\rightarrow\infty,
\end{equation}
where, here and throughout the paper, we write $a(x)\sim b(x),
\;x\to\infty$ to denote
\begin{eqnarray*}
\mathop {\lim }\limits_{x \to \infty } \frac{a(x)}{b(x)} = 1.
\end{eqnarray*}
\mbox{}\quad Many papers have been devoted to investigating the
convergence rate in $(\ref{sta})$; See Omey and Willekens
\cite{omey1986}\cite{omey1987}, Omey \cite{omey1994}, Baltr\=unas
and Omey \cite{baltrunas1998}\cite{baltrunas2002}, Baltr\=unas et
al. \cite{baltrunas2006}, Geluk and Pakes \cite{Geluk1991}, and
Geluk \cite{Geluk1992}\cite{Geluk1996}, among others. In these
papers, the precise convergence rate as well as the O-type results
has been considered. Generally speaking, the results about the
convergence rate in (\ref{sta}) would be different according to
whether or not the distribution $F$ has a finite mean.
\\ \mbox{}\quad Denote the mean of $F$ by $\mu$. In this paper, we assume $\mu<\infty$ and focus on the precise convergence rate in $(\ref{sta})$. Most of the related results
usually assume the existence of the density of $F$. For example, a
result from Omey \cite{omey1987} requires $F$ to have a
subexponential density $f$. By definition, the density $f$ is said
to be a subexponential density, denoted by $f\in \cal{S}{\mit d}$,
if
 \begin{equation}\label{lwdy} \mathop {\lim }\limits_{x \to
\infty } \frac{f(x + y)}{f(x)}=1,\quad \forall\; y\in \mathbb{R},
\end{equation}
and
\begin{eqnarray}\lim \limits_{x \to \infty }\frac{ {\int_0^x {f(y)f(x - y)dy}}}{f(x)}=
2.
\end{eqnarray}
The first part of Theorem 2.2(ii) in Omey and Willekens
\cite{omey1987} is stated as follows.
\begin{Theorem}(Omey and Willekens \cite{omey1987})\label{omthm} Suppose $E(z^{N})$ is analytic at
$z=1$, $f\in \cal{S}{\mit d}$, and
\begin{eqnarray}\overline{F^{2*}}(x)-2\overline{F}(x)\sim 2\mu
f(x),\quad x\to\infty,
\end{eqnarray}
then
\begin{equation}
\overline G (x) - (\sum\limits_{n = 0}^\infty  {np_n } )\overline F
(x) \sim \{\mu \sum\limits_{n = 0}^\infty  {n(n - 1)p_n } \}f(x),
\quad x\rightarrow \infty.
\end{equation}
\end{Theorem}
\mbox{}\quad Efforts have been taken by Omey \cite{omey1994} to
remove the condition of densities in Theorem \ref{omthm}; See
Theorem 6.1 of Omey \cite{omey1994}.  However the condition imposed
there requires $F$ to belong to a subclass of the distributions with
both dominatedly varying tails and long tails (see Omey
\cite{omey1994} for details).
\\ \mbox{}\quad In this paper, we aim to generalize Theorem \ref{omthm} to the case where the density of $F$ does not necessarily exist.
One main result of ours (see Theorem \ref{fherjdl}) unifies Theorem
\ref{omthm} and the related result in Omey \cite{omey1994}. The
appropriate condition for our result is expressed in terms of some
class of distributions, which we call the second order
subexponential distribution class. Its definition and properties
  are also stated in section 2 as main results.
The proofs are given in section 3.

\se{Main results}\label{somr} \mbox{}\quad
 Let $t\in(0,\infty]$ and write $\Delta(t) = (0,t]$,
\begin{eqnarray*}
 x + \Delta(t)  = (x,x + t]
 \end{eqnarray*}
and
\begin{eqnarray*}
 F(x + \Delta(t))=F(x,x+t]= F(x + t)-F(x).
 \end{eqnarray*}
The so-called local subexponential class as well as the local
long-tailed class is introduced by Asmussen et al.
\cite{asmussen2003}. By definition,  a distribution $F$ on $[0,
\infty )$ is said to belong to the local long-tailed class
$\mathscr{L}_{\Delta(t)}$, if the relation
\begin{equation}\label{Ldys}
F(x + y + \Delta(t) ) \sim F(x + \Delta(t) ),\quad
x\rightarrow\infty
\end{equation}
holds uniformly in $y\in [0,1]$ and hence, it holds uniformly on any
finite interval of $y$.   Furthermore, $F$ is said to belong to the
local subexponential class $\mathscr{S}_{\Delta(t)}$, if
$F\in\mathscr{L}_{\Delta(t)}$ and
\begin{equation}
F^{*2} (x + \Delta(t) ) \sim 2F(x + \Delta(t) ),\quad
x\rightarrow\infty.
\end{equation}
\begin{Definition}\label{erjczsfbdy} We say a distribution $F$ on $[0,\infty)$ with finite mean $\mu$ belongs to the second order subexponential class $\mathscr{S}_2$, if for all $t\in (0,\infty)$, $F\in\mathscr{S}_{\Delta(t)}
$ and
\begin{eqnarray}\label{ercjjgx}
\overline {F^{2*} } (x) - 2\overline F (x) \sim 2\mu F(x,x + 1],
\quad x\rightarrow \infty.
\end{eqnarray}
\end{Definition}
\begin{Proposition}\label{ncerjdl} (1) Assume $F\in\mathscr{S}_2$, then for all $n\ge 2$,
\begin{equation}\label{ncerjgx}
\overline {F^{n*} } (x) - n\overline F (x) \sim n(n - 1)\mu F(x,x +
1], \quad x\rightarrow \infty.
\end{equation}
(2) Assume for $F\in \mathscr{L}_{\Delta(t)}$ all $t\in (0,\infty)$,
$\mu<\infty$ and $ \overline F ^2 (x) = o(F(x,x + 1])$ (it means
that $\lim\limits_{x\to \infty}{\overline{F}^2(x)}/F(x,x+1]=0$). If
for some $n\ge 2$, the relation (\ref{ncerjgx}) holds, then
$F\in\mathscr{S}_2$.
\end{Proposition}
 \mbox{}\quad An uniform bound for (\ref{ncerjgx}) is given as
follows.
\begin{Proposition}\label{erjkzdl} Assume $F\in\mathscr{S}_2$, then for every fixed $\varepsilon>0$, there exist constants $A$, $K>0$, which are independent of $n$, such that for all $n\ge 2$,
\begin{equation}
\mathop {\sup }\limits_{x \ge A} \left| {\frac{{\overline {F^{n*} }
(x) - n\overline F (x)}}{{F(x,x + 1]}}} \right| \le K(1 +
\varepsilon )^n .
\end{equation}
\end{Proposition}
 \mbox{}\quad Our next result investigates the second order tail
behaviour of $G$.
\begin{Theorem}\label{fherjdl}
 (1) If $F\in\mathscr{S}_2$ and $E(z^N)$ is analytic at $z=1$, then
\begin{equation}\label{erjfhgx}
\overline G (x) - (\sum\limits_{n = 0}^\infty  {np_n } )\overline F
(x) \sim \{\mu \sum\limits_{n = 0}^\infty  {n(n - 1)p_n } \}F(x,x +
1], \quad x\rightarrow \infty.
\end{equation}
(2) Suppose $F\in \mathscr{L}_{\Delta(t)}$ for all $t\in
(0,\infty)$, $\mu<\infty$ and $ \overline F ^2 (x) = o(F(x,x + 1])$.
If the relation (\ref{erjfhgx}) holds and there exists some $l\ge 2$
such that $p_l>0$, then $F\in\mathscr{S}_2$.
\end{Theorem}
\begin{Remark} As has been shown by Asmussen et al.
\cite{asmussen2003}, if $F$ has a density $f\in \cal{S}{\mit d}$,
then for all $t\in (0,\infty)$, $F\in \mathscr{S}_{\Delta(t)}$.
Hence Theorem \ref{fherjdl} improves Theorem 2.2(ii) of Omey and
Willekens \cite{omey1987}. By Corollary \ref{condition2} (see
below), we know that Theorem \ref{fherjdl} also improves Theorem 6.1
of Omey \cite{omey1994} in the case $\mu<\infty$.
\end{Remark}
 \mbox{}\quad Next we present a result on tail equivalences.
\begin{Proposition}\label{erjwdjdl} Let $F$ and $H$ be two distributions.
If $F\in\mathscr{S}_2$ and there exist constants $K>0$, $c\in
\mathbb{R}$ such that
\begin{equation}\label{erjwdj}
\frac{\overline H (x) - K\overline F (x)}{F(x,x + 1]}\rightarrow
c,\quad x\rightarrow\infty,
\end{equation}
then $H\in\mathscr{S}_2$.
\end{Proposition}
\begin{Remark}\label{sdcz} From Proposition \ref{erjwdjdl}, we know that (\ref{erjwdj}) defines a class of distribution that is equivalent to $F$. In this equivalent class, there must exist a distribution that satisfies (\ref{erjwdj}) and has a subexponential density. To see this, let $K=1/\int_0^{1}\overline{F}(s)ds$ and define
\begin{equation}\label{hdy} \widetilde{h}(x):=KF(x,x+1],\quad \forall x>0.
\end{equation}
Assume $F\in\mathscr{S}_{\Delta(t)}$ for all $t\in (0,\infty)$. In view of Lemma \ref{lczstz} below, we have  $\widetilde{h}\in \cal{S}{\mit d}$. Denote the distribution function of $\widetilde{h}(x)$ by $H$. It is easy to see
\begin{eqnarray}\int_x^{\infty}F(s,s+1]ds=\int_x^{x+1}\overline{F}(s)ds=\int_0^1\overline{F}(x+z)dz.
\end{eqnarray}
Since $F\in \mathscr{L}_{\Delta(t)}$ for all $t\in
(0,\infty)$, then it follows from Lemma \ref{gtdlyl} below and the dominated convergence theorem that \begin{eqnarray}\overline{H}(x)-K\overline{F}(x)&=&-K\int_0^1(\overline{F}(x)-\overline{F}(x+z))dz\nonumber\\&\sim &-K F(x,x+1]\int_0^1zdz\nonumber\\&=&-\frac{K}{2}F(x,x+1],
\end{eqnarray}
i.e., $H$ satisfies (\ref{erjwdj}) with $c=-K/2$.
\end{Remark}
\begin{Remark} It follows from Proposition \ref{ncerjdl}(1) that in Proposition \ref{erjwdjdl}, (\ref{erjwdj}) implies that (\ref{erjwdj}) holds with $H$ and $F$ (in the numerator) replaced by $H^{n*}$ and $F^{n*}$ and $c$ replaced by $nc+K(\mu_H-\mu)n(n-1)$, where $\mu_H=\int_0^{\infty}\overline{H}(x)dx<\infty$.
\end{Remark}
\mbox{}\quad The following lemma about local subexponential distributions, which is cited by Remark \ref{sdcz}, might be of independent interest.
\begin{Lemma}\label{lczstz} Let $t\in (0,\infty)$ be fixed, then $F\in \mathscr{S}_{\Delta(t)}$ if and only if
$KF(\cdot+\Delta(t))\in \cal{S}{\mit d}$, where $K$ as a positive constant, is defined as \begin{eqnarray}K=\frac{1}{\int_0^{t}\overline{F}(s)ds}.
\end{eqnarray}
\end{Lemma}
 \mbox{}\quad Finally, we give some sufficient conditions for $F\in
\mathscr{S}_2$. A distribution $F$ on $[0,\infty)$ is said to belong
to $\mathscr{S}^{*}$ (see Kl\"{u}ppelberg \cite{kluppelberg1988}),
if
\begin{eqnarray}
\int_{0}^{x}\overline{F}(y)\overline{F}(x-y)dy\sim
2\mu\overline{F}(x),\quad x\rightarrow\infty.
\end{eqnarray}
It is well known that $\mathscr{S}^{*}\subset \mathscr{S}$. Denote
$h(x)=F(x,x+1]$ and $q(x)=h(x)/\overline{F}(x)$.
\begin{Proposition}\label{condition1} Suppose $F\in
\mathscr{L}_{\Delta(t)}$ for all $t\in (0,\infty)$, $\mu<\infty$,
$F\in \mathscr{S}^{*}$, $ \overline F ^2 (x/2) = o(F(x,x + 1])$ and
for all $y>0$,
\begin{eqnarray}\label{qd}\limsup\limits_{x\to\infty}\frac{q(xy)}{q(x)}<\infty.
\end{eqnarray}
 Then $F\in\mathscr{S}_2$.
\end{Proposition}
\begin{Remark} In view of Proposition \ref{ncerjdl}(1), we know that
Proposition \ref{condition1} improves Proposition 3.5(iii) of
Baltr\={u}nas \cite{baltrunas1999}.
\end{Remark}
\begin{Corollary}\label{condition2} Suppose $F\in
\mathscr{L}_{\Delta(t)}$ for all $t\in (0,\infty)$, $\mu<\infty$, $
\overline F ^2 (x/2) = o(F(x,x + 1])$ and for all $y>0$,
\begin{eqnarray}\label{hd}\limsup\limits_{x\to\infty}\frac{h(xy)}{h(x)}<\infty.
\end{eqnarray}
 Then $F\in\mathscr{S}_2$.
\end{Corollary}
\mbox{}\quad Some typical subexponential distributions including the Pareto, lognormal and Weibull (with parameter between 0 and 1) distributions all belong to $\mathscr{S}_2$, which is shown in the following.\\
\mbox{}\quad For the Pareto distribution $F$, i.e.,
 $\overline{F}(x)=cx^{-\alpha}$, where $c>0$ and $\alpha>1$, it is easy to obtain for every fixed $t\in (0,\infty)$,
 \begin{eqnarray}\label{prf} F(x,x+t]\sim c\alpha t x^{-(\alpha+1)},\quad x\to\infty,
 \end{eqnarray}
and hence by Corollary \ref{condition2}, it is easy to see $F\in\mathscr{S}_2$.\\
\mbox{}\quad Let $F$ be the lognormal distribution with the density $f(x)=e^{-(lnx-\mu)^2/2\sigma^2}/x\sqrt {2\pi \sigma ^2 }$.
Let $\Phi$ be the standard normal distribution with the density $\phi$. Then by using the relation between the lognormal and normal distributions, and the following well-known relation
\begin{eqnarray}1-{\Phi}(x)\sim\frac{1}{x}\phi(x),\quad x\to\infty,
\end{eqnarray}
it is easy to obtain
\begin{eqnarray}\overline{F}(x)=1-{\Phi}(\frac{lnx-\mu}{\sigma})\sim \frac{\sigma}{lnx}{\phi}(\frac{lnx-\mu}{\sigma}),\quad x\to\infty.
\end{eqnarray}
On the other hand, it is easy to see for every fixed $t\in (0,\infty)$,
\begin{eqnarray}F(x,x+t]\sim t f(x)=\frac{t}{x\sigma}\phi(\frac{lnx-\mu}{\sigma}),\quad x\to\infty.
\end{eqnarray}
Thus,
\begin{eqnarray}q(x)\sim \frac{lnx}{x}\cdot \frac{1}{\sigma^2}, \quad x\to\infty.
\end{eqnarray}
By Proposition \ref{condition1}, it is easy to see $F\in\mathscr{S}_2$.\\
\mbox{}\quad For the Weibull distribution $F$, i.e., $\overline{F}(x)=e^{-x^{\beta}}$, $\beta\in (0,1)$, we have for every fixed $t\in (0,\infty)$,
\begin{eqnarray}F(x,x+t]=\beta t x^{\beta-1}e^{-x^{\beta}}, \quad x\to\infty.
\end{eqnarray}
Hence
\begin{eqnarray}q(x)\sim \beta x^{\beta-1}, \quad x\to\infty.
\end{eqnarray}
By Proposition \ref{condition1}, it is easy to see $F\in\mathscr{S}_2$.\\
\mbox{}\quad A distribution, which belongs to $\mathscr{S}_2$ but does not have a density, is presented in the following example.
\begin{Example} Define for $n\ge 2$,
\begin{eqnarray}\overline{F}(x)=c(1+\frac{1}{n})x^{-\alpha}, \quad n^{\beta}\le x <(n+1)^{\beta},
\end{eqnarray}
where $c>0$, $\alpha>1$ and $\beta\in (1,2)$. Since
\begin{eqnarray} (n+1)^{\beta}-n^{\beta}=n^{\beta}[(1+\frac{1}{n})^{\beta}-1]\sim\beta n^{\beta-1}\to \infty,\quad n\to\infty,
\end{eqnarray}
then for any fixed $t\in (0,\infty)$ and sufficiently large $x$, there only exist two cases: $n^{\beta}\le x<x+t<(n+1)^{\beta}$ or $n^{\beta}\le x<(n+1)^{\beta}\le x+t<(n+2)^{\beta}$ for some $n$. In either case, through some simple calculations, it is to easy see that the relation (\ref{prf}) always holds. From this and in view of
\begin{eqnarray}\overline{F}(x)\sim cx^{-\alpha},\quad x\to\infty,
\end{eqnarray}
it is easy to see that the conditions of Corollary \ref{condition2} are satisfied, and thus,
 $F\in\mathscr{S}_2$.
 However, since $F$ is not continuous,  it does not have a density.
\end{Example}

\se{Proofs}\label{proof} \mbox{}\quad In the sequel, all limit
relations between two functions $g_1(x)$ and $g_2(x)$ of one
variable $x$, unless explicitly stated otherwise, are for
$x\to\infty$.
 If  $g_1$ or $g_2$ is a function
of two variables $x$ and $A$, then the limit relations between them,
unless explicitly stated otherwise, are for $x\rightarrow\infty$ and
then $A\rightarrow\infty$, the meaning of which is specified as
follows:\\
 $g_1=o(g_2)$ denotes $$ \mathop {\lim }\limits_{A \to \infty }
\mathop {\lim \sup }\limits_{x \to \infty } \left| g_1 /g_2 \right|
= 0; $$\\  $g_1\sim g_2$ denotes $$ \mathop {\lim }\limits_{A \to
\infty } \mathop {\lim \sup }\limits_{x \to \infty } \left| {g_1
/g_2 - 1} \right| = 0; $$ \\
 $g_1\lesssim g_2$ denotes $$ \mathop {\lim \sup }\limits_{A \to
\infty } \mathop {\lim \sup }\limits_{x \to \infty } g_1 /g_2 <
\infty; $$\\  $g_1\gtrsim  g_2$ denotes $$ \mathop {\lim \inf
}\limits_{A \to \infty } \mathop {\lim \inf }\limits_{x \to \infty }
g_1 /g_2 >0. $$
\begin{Lemma}\label{gtdlyl} Assume $F\in
\mathscr{L}_{\Delta(t)}$ for all $t\in (0,\infty)$. Then for all
$t\in (0,\infty)$,
\begin{equation}\label{gtdlyl2}
\frac{{F(x + \Delta(t) )}}{{F(x,x + 1]}} \to t,\begin{array}{*{20}c}
   {} & {x \to \infty .}  \\
\end{array}
\end{equation}
\end{Lemma}
{\bf Proof.} For any $ \delta  \in (0,\min \{ t,1\} )$, there exist
positive integers $k$, $n$ such that
\begin{eqnarray}
k \delta  \le t < (k+1) \delta ,\begin{array}{*{20}c}
   {} & {}  \\
\end{array}n \delta  \le 1 < (n+1) \delta .
\end{eqnarray}
Obviously, when $\delta\rightarrow 0+$,
\begin{eqnarray}\label{knjj}
k \sim \frac{t}{\delta },\begin{array}{*{20}c}
   {} & {}  \\
\end{array}n \sim \frac{1}{\delta }.
\end{eqnarray}
Obviously,
\begin{eqnarray}
&&\sum\limits_{i = 1}^k {F(x + (i - 1)\delta ,x + i\delta ]}  \le
F(x + \Delta(t) ) \le \sum\limits_{i = 1}^{k + 1} {F(x + (i -
1)\delta ,x +
i\delta ]},\nonumber\\
&&\sum\limits_{i = 1}^n {F(x + (i - 1)\delta ,x + i\delta ]}  \le
F(x,x + 1] \le \sum\limits_{i = 1}^{n + 1} {F(x + (i - 1)\delta ,x +
i\delta ]}.
\end{eqnarray}
Let $\delta$ be fixed, then for all $i=1$, 2, $\cdots$, $\max\{k,
n\}$,
\begin{eqnarray}F(x + (i - 1)\delta ,x +
i\delta ]\sim F(x ,x + \delta ],
\end{eqnarray}
and hence
\begin{eqnarray}\label{1a1a}
\frac{k}{{n + 1}} \le \mathop {\lim \inf }\limits_{x \to \infty }
\frac{{F(x + \Delta(t) )}}{{F(x,x + 1]}} \le \mathop {\lim \sup
}\limits_{x \to \infty } \frac{{F(x + \Delta(t) )}}{{F(x,x + 1]}}
\le \frac{{k + 1}}{n}.
\end{eqnarray}
Let $\delta\rightarrow 0+$ in (\ref{1a1a}) and in view of
(\ref{knjj}), we obtain (\ref{gtdlyl2}).
\begin{Lemma}\label{jbcfby} For any $t\in (0,\infty)$, the following three assertions are equivalent:\\
(1) $F\in \mathscr{S}_{\Delta(t)}$,\\
(2) $F\in \mathscr{L}_{\Delta(t)}$ and
\begin{eqnarray}\label{djs1}
\int_0^{x - A} {F(x - y + \Delta(t) )d} F(y) \sim F(x + \Delta(t) ),
\end{eqnarray}
(3) $F\in \mathscr{L}_{\Delta(t)}$ and
\begin{eqnarray}
\int_A^{x - A} {F(x - y + \Delta(t) )d} F(y) =o( F(x + \Delta(t) )).
\end{eqnarray}
\end{Lemma}
{\bf Proof.} The proof of this lemma is similar to that of
Proposition 2 of Asmussen et al. \cite{asmussen2003}, so we omit it.
\begin{Lemma}\label{erjdjyl} Assume $F\in \mathscr{L}_{\Delta(t)}$ for all $t\in
(0,\infty)$ and $\mu<\infty$. Then the relation (\ref{ercjjgx}) is
equivalent to
\begin{equation}
\int_A^{x - A} {\{\overline F (x - y) - \overline F (x)\}} dF(y) -
\overline F ^2 (x) = o(F(x,x + 1]).
\end{equation}
\end{Lemma}
{\bf Proof.} Assume $F\in \mathscr{L}_{\Delta(t)}$ for all $t\in
(0,\infty)$ and $\mu<\infty$. Notice that
\begin{eqnarray}\label{2cfjs}
\overline {F^{2*} } (x) - 2\overline F(x) = \int_0^x {\{\overline F
(x - y) - \overline F (x)\}} dF(y) - \overline F^2 (x).
\end{eqnarray}
By Lemma \ref{gtdlyl}, it is obvious that
\begin{eqnarray}\label{2cfjcf1}
\int_0^A {\{\overline F (x - y) - \overline F (x)\}} dF(y) \sim
\int_0^{\infty} y dF(y)\cdot F(x,x + 1]=\mu F(x,x + 1].
\end{eqnarray}
By integrating by parts, we obtain
\begin{eqnarray}\label{2cfjcf3}
&&\int_{x - A}^x {\{\overline F (x - y) - \overline F (x)\}} dF(y)
\nonumber\\
&=& \int_0^A {\{\overline F (x - y) - \overline F (x)\}} dF(y) +
\{\overline F (x - A) - \overline F (x)\}\{\overline F (A) -
\overline F (x)\},
\end{eqnarray}
hence by Lemma \ref{gtdlyl} and in view of the fact that $\mathop
{\lim }\limits_{A \to \infty } A\overline F (A) = 0$, we have
\begin{eqnarray}\label{2cfjcf2}
\int_{x - A}^x {\{\overline F (x - y) - \overline F (x)\}} dF(y)
\sim \mu F(x,x + 1].
\end{eqnarray}
By (\ref{2cfjs}), (\ref{2cfjcf1}) and (\ref{2cfjcf2}), we obtain the
desired result.
\begin{Lemma}\label{T1lxtl} Assume $F\in
\mathscr{L}_{\Delta(t)}$  for all $t\in (0,\infty)$ and $\overline F
^2 (x) = o(F(x,x + 1])$. Then the relation (\ref{ercjjgx}) implies
$F\in \mathscr{S}_{\Delta(t)}$ for all $t\in (0,\infty)$.
\end{Lemma}
{\bf Proof.} Assume the relation (\ref{ercjjgx}) holds. Since
$\overline{F} ^2 (x)= o(F(x,x+1])$, from Lemma \ref{erjdjyl} it
follows that
\begin{eqnarray}\label{er2}\int_A^{x-A} {\{\overline F (x - y) - \overline F
(x)\}}dF(y)=o(F(x,x+1]).
\end{eqnarray}
Hence for $x>2A$ and $A>t$,
\begin{eqnarray}
\int_A^{x - A} {F(x - y + \Delta(t) )d} F(y)&\le & \int_A^{x-A}
{\{\overline F (x - y) - \overline F
(x)\}}dF(y)\nonumber\\&=&o(F(x,x+1]) =o( F(x + \Delta(t) )).
\end{eqnarray}
Thus by lemma \ref{jbcfby}, we prove $F\in \mathscr{S}_{\Delta(t)}$.
\\
{\bf Proof of Proposition \ref{ncerjdl}(1).}  We argue by induction.
First the relation (\ref{ncerjgx}) is trivial for $n=2$.
Furthermore, assume (\ref{ncerjgx}) holds for some $n-1\ge 2$, i.e.,
\begin{equation}\label{n-1cerjgx}
\overline {F^{(n-1)*} } (x) - (n-1)\overline F (x) \sim (n-1)(n -
2)\mu F(x,x + 1].
\end{equation}
Then it suffices to prove (\ref{ncerjgx}) for $n$.  Note that
\begin{eqnarray}\label{ncerjfj}
&&\overline {F^{n*} } (x) - n\overline F(x)\nonumber\\
&=& \int_0^x {\{\overline {F^{(n - 1)*} } (x - y) - (n - 1)\overline
F(x
- y)\}} dF(y)\nonumber\\&& + (n - 1)\{\overline {F^{2*} } (x) - 2\overline F(x)\}\nonumber\\
&:=& I_1  + I_2 .
\end{eqnarray}
Obviously,
\begin{equation}
I_2  \sim 2(n - 1)\mu F(x,x + 1].
\end{equation}
For $x>A>0$,
\begin{eqnarray}
I_1  &= &\int_0^{x - A} {\{\overline {F^{(n - 1)*} } (x - y) - (n -
1)\overline F(x - y)\}} dF(y)\nonumber\\
&& + \int_{x - A}^x {\{\overline {F^{(n - 1)*} } (x - y) - (n -
1)\overline F(x - y)\}} dF(y)\nonumber\\&: =& J_1  + J_2 .
\end{eqnarray}
Since $F\in \mathscr{S}_{\Delta(1)}$, it follows from Lemma
\ref{jbcfby} that the relation (\ref{djs1}) holds for $t=1$. Hence
by (\ref{n-1cerjgx}), we obtain
\begin{eqnarray}
J_1  &\sim & (n - 1)(n - 2)\mu \int_0^{x - A} {F(x - y,x - y + 1]}
dF(y)\nonumber\\
& \sim & (n - 1)(n - 2)\mu F(x,x + 1].
\end{eqnarray}
For $J_2$, by integrating by parts, we obtain
\begin{eqnarray}
J_2  &=& \int_0^A {\{\overline F (x - y) - \overline F (x)\}} dF^{(n
- 1)*} (y)- (n - 1)\int_0^A {\{\overline F (x - y) - \overline F
(x)\}}
dF(y)\nonumber\\
&& + \{\overline {F^{(n - 1)*} } (A) - (n - 1)\overline
F(A)\}\{\overline
F(x - A) - \overline F (x)\}\nonumber\\
&: =& K_1  - K_2  + K_3 .
\end{eqnarray}
By Lemma \ref{gtdlyl}, it follows that
\begin{eqnarray}
K_1  \sim \int_0^\infty  y dF^{(n - 1)*} (y)\cdot F(x,x + 1] = (n -
1)\mu F(x,x + 1],
\end{eqnarray}
\begin{eqnarray}
K_2  \sim (n - 1)\int_0^\infty  y dF(y)\cdot F(x,x + 1] = (n - 1)\mu
F(x,x + 1],
\end{eqnarray}
and
\begin{eqnarray}\label{ncerjfjzhys}
K_3  \sim \{\overline {F^{(n - 1)*} } (A) - (n - 1)\overline F(A)\}A
\cdot F(x,x + 1] = o(F(x,x + 1]).
\end{eqnarray}
Then we have
\begin{eqnarray}\label{j2o}J_2=o(F(x,x+1]),
\end{eqnarray}
and hence
\begin{eqnarray}
\overline {F^{n*} } (x) - n\overline F(x) &\sim & (n - 1)(n - 2)\mu
F(x,x
+ 1]\nonumber\\
&& + 2(n - 1)\mu F(x,x + 1] = n(n - 1)\mu F(x,x + 1],
\end{eqnarray}
as required.
\\
 \mbox{}\quad The proof of Proposition \ref{ncerjdl}(2) needs the
following Lemma.
\begin{Lemma}\label{erjjxgxyl} Assume $F\in \mathscr{L}_{\Delta(t)}$ for all $t\in (0,\infty)$,  $\mu<\infty$ and $
\overline F ^2 (x) = o(F(x,x + 1])$, then for all $n\ge 2$,
\begin{eqnarray}\label{erjjxgx}
\mathop {\lim \inf }\limits_{x \to \infty } \frac{{\overline {F^{n*}
} (x) - n\overline{F}(x)}}{{F(x,x + 1]}} \ge n(n - 1)\mu .
\end{eqnarray}
\end{Lemma}
{\bf Proof.} We still argue by induction. In the following, we use the same notations ($J_1$, $I_2$, $\cdots$) as in the proof of Proposition \ref{ncerjdl}(1). From (\ref{2cfjcf1}) and (\ref{2cfjcf2}), it follows that for $x>2A$,
\begin{eqnarray}
&&\int_0^x {\{\overline F (x - y) - \overline F
(x)\}dF(y)}\nonumber\\
& \ge & \int_0^A {\{\overline F (x - y) - \overline F (x)\}dF(y)}
+\int_{x-A}^x {\{\overline F (x - y) - \overline F (x)\}dF(y)}
\nonumber\\
 &\sim & 2\mu F(x,x + 1].
\end{eqnarray}
Then, in view of (\ref{2cfjs}) and the condition $\overline F ^2 (x)
= o(F(x,x + 1])$, we prove (\ref{erjjxgx}) for $n=2$. Assume
(\ref{erjjxgx}) holds for some $n-1\ge2$, i.e.,
\begin{eqnarray}\label{n-1erjjxgx}
\mathop {\lim \inf }\limits_{x \to \infty } \frac{{\overline
{F^{(n-1)*} } (x) - (n-1)\overline{F}(x)}}{{F(x,x + 1]}} \ge (n-1)(n
- 2)\mu .
\end{eqnarray}
Then for $x>2A$, we
 have
\begin{eqnarray}\label{ii1jx}
J_1  &\gtrsim & (n - 1)(n - 2)\mu \int_0^{x - A} {F(x - y,x - y +
1]}
dF(y)\nonumber\\
&\ge & (n - 1)(n - 2)\mu \int_0^{A} {F(x - y,x - y + 1]}dF(y)\nonumber\\
 &\sim & (n - 1)(n - 2)\mu F(x,x + 1].
\end{eqnarray}
From the proof of Proposition \ref{ncerjdl}(1), we have (\ref{j2o}).
Moreover, the relation (\ref{erjjxgx}) holds for $n=2$, i.e.
\begin{equation}\label{i2jx}
I_2  = (n - 1)\{\overline {F^{2*} } (x) - 2\overline F(x)\} \gtrsim
2(n - 1)\mu F(x,x + 1].
\end{equation}
Hence we have
\begin{eqnarray}
\overline {F^{n*} } (x) - n\overline F(x) &\gtrsim  & (n - 1)(n -
2)\mu F(x,x
+ 1]\nonumber\\
&& + 2(n - 1)\mu F(x,x + 1] = n(n - 1)\mu F(x,x + 1].
\end{eqnarray}
{\bf Proof of Proposition \ref{ncerjdl}(2).} If the relation
(\ref{ncerjgx}) holds for $n=2$, the result is obvious. Thus, we
assume the relation (\ref{ncerjgx}) holds for some $n\ge 3$. From
the proof of Lemma \ref{erjjxgxyl}, we know that the relations from
(\ref{n-1erjjxgx}) to (\ref{i2jx}) still hold. However the relation
(\ref{ncerjgx}) implies
\begin{eqnarray}J_1+J_2+I_2\sim 2n(n - 1)\mu F(x,x + 1],
\end{eqnarray}
hence (\ref{ii1jx}) and (\ref{i2jx}) necessarily hold with the sign
$\gtrsim$ replaced by $\sim$. In particular, we have
\begin{eqnarray}\label{I2sim} I_2\sim 2(n
- 1)\mu F(x,x + 1],
\end{eqnarray}
which is equivalent to (\ref{ercjjgx}). From this and Lemma
\ref{T1lxtl}, we have $F\in \mathscr{S}_{\Delta(t)}$ for all $t\in
(0,\infty)$ and hence, the proof is completed.
\\
{\bf Proof of Proposition \ref{erjkzdl}.} Without loss of
generality, we assume $\varepsilon \in(0,1)$. By Lemma
 \ref{jbcfby}, we know there exist sufficiently large constants $A$, $A^{'}$ such that
 $A>A^{'}>0$ and
\begin{eqnarray}\label{jbkzbds2}
\mathop {\sup }\limits_{x \ge A}\left\{{{\int_0^{x - A^{'}} {F(x -
y,x - y + 1]} dF(y)} \mathord{\left/
 {\vphantom {{\int_0^{x - A} {F(x - y,x - y + 1]} dF(y)} {F(x,x + 1]}}} \right.
 \kern-\nulldelimiterspace} {F(x,x + 1]}}\right\} \le 1 +
 \varepsilon/4
\end{eqnarray}
and
\begin{eqnarray}\label{erjkzgx1}
\mathop {\sup }\limits_{x \ge A} \left\{ {{{\left|\overline {F^{2*}
} (x) - 2\overline F(x)\right|} \mathord{\left/
 {\vphantom {{\{\overline {F^{2*} } (x) - 2\overline F(x)\}} {F(x,x + 1]}}} \right.
 \kern-\nulldelimiterspace} {F(x,x + 1]}}} \right\} < 3\mu.
\end{eqnarray}
Obviously,
\begin{eqnarray}\int_0^{x - A} {F(x -
y,x - y + 1]} dF(y)\le\int_0^{x - A^{'}} {F(x - y,x - y + 1]} dF(y),
\end{eqnarray}
Hence by (\ref{jbkzbds2}), we know that
\begin{eqnarray}\label{jbkzbds}\mathop {\sup }\limits_{x \ge A}\left\{{{\int_0^{x - A} {F(x -
y,x - y + 1]} dF(y)} \mathord{\left/
 {\vphantom {{\int_0^{x - A} {F(x - y,x - y + 1]} dF(y)} {F(x,x + 1]}}} \right.
 \kern-\nulldelimiterspace} {F(x,x + 1]}}\right\} \le 1 +
 \varepsilon/4.
\end{eqnarray}
Since $F(\log{x},\log{x} + 1]$ is a slowly varying function, so is
$1/F(\log{x},\log{x} + 1]$, hence by Lemma 1.3.2 of Bingham et al.
\cite{bingham1987}, the above $A$ can be chosen such that the
function $1/F(\log{x},\log{x} + 1]$ is locally bounded on $[e^A,
\infty)$, i.e., $1/F(x,x + 1]$ is locally bounded on $[A, \infty)$.
Hence by Lemma
 \ref{gtdlyl}, we know there exists a sufficiently large constant
 $B>A$ such that
\begin{eqnarray}\label{kz1}
\mathop {\sup }\limits_{x \ge B} \left\{ {{{\{\overline F (x - A) -
\overline F (x)\}} \mathord{\left/
 {\vphantom {{[\overline F (x - A) - \overline F (x)]} {F(x,x + 1]}}} \right.
 \kern-\nulldelimiterspace} {F(x,x + 1]}}} \right\}<\infty,
\end{eqnarray}
and
\begin{eqnarray}\label{kz2}
\mathop {\sup }\limits_{A \le x < B} \left\{ {{{\{\overline F (x -
A) - \overline F (x)\}} \mathord{\left/
 {\vphantom {{\{\overline F (x - A) - \overline F (x)\}} {F(x,x + 1]}}} \right.
 \kern-\nulldelimiterspace} {F(x,x + 1]}}} \right\} \le \mathop {\sup }\limits_{A \le x < B} \left\{ {{1 \mathord{\left/
 {\vphantom {1 {F(x,x + 1]}}} \right.
 \kern-\nulldelimiterspace} {F(x,x + 1]}}} \right\} < \infty .
\end{eqnarray}
Thus, there exists a positive constant $M$, which is independent of
$n$, such that both the left-hand sides of (\ref{kz1}) and
(\ref{kz2}) do not exceed $M$. On the other hand, by the definition
of $J_2$, it is easy to see
\begin{eqnarray}
\left|J_2\right|\le n \left\{ \overline{F} (x - A) - \overline{F}(x)
\right\}.
\end{eqnarray}
Hence we have
\begin{eqnarray}\label{erjkzgx2}
\mathop {\sup }\limits_{x \ge A} \left\{ {{{\left|J_2\right| }
\mathord{\left/
 {\vphantom {{J_2 } {F(x,x + 1]}}} \right.
 \kern-\nulldelimiterspace} {F(x,x + 1]}}} \right\} &\le & n \mathop {\sup }\limits_{x \ge A} \left\{ {{{\{\overline F (x - A) - \overline F (x)\}} \mathord{\left/
 {\vphantom {{\{\overline F (x - A) - \overline F (x)\}} {F(x,x + 1]}}} \right.
 \kern-\nulldelimiterspace} {F(x,x + 1]}}} \right\}\nonumber\\
 &\le & Mn < \infty .
\end{eqnarray}
 Denote
\begin{eqnarray} \alpha _n  = \mathop {\sup }\limits_{x \ge A
}\left|\frac{{\overline {F^{n*} } (x) - n\overline F(x)}}{{F(x,x +
1]}}\right|.
\end{eqnarray}
By (\ref{jbkzbds}), we have
\begin{eqnarray}\label{erjkzgx3}
&&\mathop {\sup }\limits_{x \ge A} \left\{ {{{\left|J_1\right| }
\mathord{\left/
 {\vphantom {{J_1 } {F(x,x + 1]}}} \right.
 \kern-\nulldelimiterspace} {F(x,x + 1]}}} \right\}\nonumber\\ &\le &\alpha _{n - 1} \mathop {\sup }\limits_{x \ge A} \left\{ {{{\int_0^{x - A} {F(x - y,x - y + 1]} dF(y)} \mathord{\left/
 {\vphantom {{\int_0^{x - A} {F(x - y,x - y + 1]} dF(y)} {F(x,x + 1]}}} \right.
 \kern-\nulldelimiterspace} {F(x,x + 1]}}} \right\}\nonumber\\
 & \le & (1 + \varepsilon/4 )\alpha _{n - 1}.
\end{eqnarray}
From (\ref{erjkzgx1}), (\ref{erjkzgx2}) and (\ref{erjkzgx3}) it
follows that
\begin{eqnarray}
\alpha _n  \le (1 + \varepsilon/4 )\alpha _{n - 1}  + 3\mu (n - 1) +
Mn \le (1 + \varepsilon/4 )\alpha _{n - 1}  + C_1n,
\end{eqnarray}
where $C_1=3\mu+M$. By induction and in view of $\alpha_1=0$, we
obtain
\begin{eqnarray}\label{akz}
\alpha _n  \le C_1\sum\limits_{i = 0}^{n - 2} {(n - i)(1 +
\varepsilon /4)^i }  \le C_1n^2 (1 + \varepsilon /4)^n ,
\end{eqnarray}
It is easy to see that the right-hand side of (\ref{akz}) does not
exceed $K(1 + \varepsilon )^n$ for an appropriately chosen constant
$K$ and hence, the proof is completed.
\\
\mbox{}\quad Let
\begin{eqnarray}
\beta _n  = \mathop {\inf }\limits_{x \ge A } \frac{{\overline
{F^{n*} } (x) - n\overline F(x)}}{{F(x,x + 1]}}.
\end{eqnarray}
\begin{Lemma}\label{erjxjkzyl} Assume $F$ is a distribution on $[0,\infty)$ satisfying $ \overline F ^2
(x) = o(F(x,x + 1])$. Then there exists a constant $A>0$, which is
independent of $n$, such that for all $n\ge 2$,
\begin{eqnarray}\label{betainf}
\beta _n  \ge - n^2.
\end{eqnarray}
\end{Lemma}
{\bf Proof.} By Bonfferoni's inequality, we have
\begin{eqnarray}\label{bfi}
\overline {F^{n*} } (x)& =& P\left( {S_n  > x} \right)\nonumber\\
&\ge & P\left( {\mathop {\max }\limits_{1 \le k \le n} X_k  > x}
\right)\nonumber\\
& \ge & \sum\limits_{k = 1}^n {P\left( {X_k  > x} \right)}  -
\sum\limits_{1 \le i < j \le n} {P\left( {X_i  > x,X_j
> x} \right)}\nonumber\\
&\ge & n\overline F (x) - n^2 \overline F ^2 (x),
\end{eqnarray}
Since $ \overline F ^2 (x) = o(F(x,x + 1])$, there exists a
sufficiently large constant $A>0$ such that
\begin{eqnarray}\label{Fkzgx}
\sup\limits_{x\ge A}{\left\{\overline{F}^2 (x)/ F(x,x +
1]\right\}}\le 1.
\end{eqnarray}
Combining (\ref{bfi}) and (\ref{Fkzgx}) gives (\ref{betainf}).
\\
{\bf Proof of Theorem \ref{fherjdl}.} (1)
By Proposition \ref{ncerjdl}(1), Proposition \ref{erjkzdl} and the dominated convergence theorem, we obtain the desired result.\\
(2) Obviously,
\begin{eqnarray}\label{plkz1}
&& p_l \mathop {\lim \sup }\limits_{x \to \infty } \frac{{\overline
{F^{l*} } (x) - l\overline F(x)}}{{F(x,x + 1]}}\nonumber\\
 &\le & \mathop
{\lim }\limits_{x \to \infty } \frac{{\overline G(x) -
(\sum\limits_{n = 0}^\infty  {np_n } )\overline F(x)}}{{F(x,x + 1]}}
- \mathop {\lim \inf }\limits_{x \to \infty } \frac{{\sum\limits_{n
\ne l} {\{\overline
{F^{n*} } (x) - n\overline F(x)\}} }p_n}{{F(x,x + 1]}} \nonumber\\
 &= & \mu \sum\limits_{n = 0}^\infty  {n(n - 1)p_n } - \mathop {\lim \inf
}\limits_{x \to \infty } \frac{{\sum\limits_{n \ne l} {\{\overline
{F^{n*} } (x) - n\overline F(x)\}} }p_n}{{F(x,x + 1]}}.
\end{eqnarray}
By Lemma \ref{erjxjkzyl}, we know that Fatou's Lemma (cf. p. 94 of
Chow and Teicher \cite{chow1978} ) can be applied to the second term
above, which gives
\begin{eqnarray}\label{plkz2}&&\mathop {\lim \inf
}\limits_{x \to \infty } \frac{{\sum\limits_{n \ne l} {\{\overline
{F^{n*} } (x) - n\overline F(x)\}} }p_n}{{F(x,x + 1]}}\nonumber\\
&\ge & \sum\limits_{n \ne l}\mathop {\lim \inf }\limits_{x \to
\infty }\left\{ \frac{{ {\overline {F^{n*} } (x) - n\overline F(x)}
}}{{F(x,x +
1]}} \right\}p_n \nonumber\\
&\ge&\mu \sum\limits_{n\neq l} {n(n - 1)p_n },
\end{eqnarray}
where in the last step, Lemma \ref{erjjxgxyl} has been applied.
Combining (\ref{plkz1}) and (\ref{plkz2}) gives
\begin{eqnarray}p_l \mathop {\lim \sup }\limits_{x \to \infty } \frac{{\overline
{F^{l*} } (x) - l\overline F(x)}}{{F(x,x + 1]}}&\le & \mu
\sum\limits_{n = 0}^\infty  {n(n - 1)p_n }- \mu \sum\limits_{n\neq
l} {n(n - 1)p_n }\nonumber\\ &=& p_l l(l - 1)\mu.
\end{eqnarray}
From this and Lemma \ref{erjjxgxyl} it follows that
\begin{eqnarray}
\mathop {\lim }\limits_{x \to \infty } \frac{{\overline {F^{l*} }
(x) - l\overline F(x)}}{{F(x,x + 1]}} = l(l - 1)\mu.
\end{eqnarray}
Hence by Proposition \ref{ncerjdl}(2), we obtain $F\in
\mathscr{S}_2$.
\\
{\bf Proof of Proposition \ref{erjwdjdl}.} Notice that
\begin{eqnarray}\label{bzxjc}
\frac{{H(x + \Delta(t) )}}{{K \cdot F(x + \Delta(t) )}} - 1 =
\frac{{\overline H (x) - K \cdot \overline F (x)}}{{K \cdot F(x +
\Delta(t) )}} - \frac{{\overline H (x + t) - K \cdot \overline F (x
+ t)}}{{K \cdot F(x + \Delta(t) )}}
\end{eqnarray}
and
\begin{eqnarray}
F(x+t+\Delta(t) )\sim F(x + \Delta(t) ).
\end{eqnarray}
Hence by (\ref{erjwdj}) and Lemma \ref{gtdlyl}, we know that the
right-hand side of (\ref{bzxjc}) tends to zero, i.e.,
\begin{eqnarray}\label{GFjbdj}
H(x + \Delta(t) ) \sim K \cdot F(x + \Delta(t) ).
\end{eqnarray}
Hence by Lemma 1 of Asmussen et al. \cite{asmussen2003}, we have
\begin{eqnarray}\label{Gjbczsx222}H\in \mathscr{S}_{\Delta(t)}.
\end{eqnarray}
By (\ref{erjwdj}),
\begin{eqnarray}\label{pfgx}
\overline H ^2 (x) - K^2  \cdot \overline F ^2 (x) =\{\overline H(x)
- K \cdot \overline F (x)\}\{\overline H (x) + K \cdot \overline F
(x)\} = o(F(x,x + 1]).
\end{eqnarray}
Notice that
\begin{eqnarray}\label{erjjhfj}
&&\int_A^{x - A} {\{\overline H (x - y) - \overline H(x)\}}
dF(y)\nonumber\\
& =& K\int_A^{x - A} {\{\overline F (x - y) - \overline F (x)\}}
dF(y) + \int_A^{x - A} {\{\overline H (x - y) - K \cdot \overline F
(x -
y)\}} dF(y)\nonumber\\
&& - \int_A^{x - A} {\{\overline H (x) - K \cdot \overline F (x)\}}
dF(y).
\end{eqnarray}
By (\ref{erjwdj}), we have
\begin{eqnarray}\label{erjjhfjcf1}
\int_A^{x - A} {\{\overline H (x) - K \cdot \overline F (x)\}}
dF(y)\le \{\overline H(x) - K \cdot \overline F (x)\}\overline{F}(A)
= o(F(x,x + 1])
\end{eqnarray}
and
\begin{eqnarray}\label{erjjhfjcf2}
&&\int_A^{x - A} {\{\overline H (x - y) - K \cdot \overline F (x -
y)\}} dF(y)\nonumber\\
 &\sim & c\int_A^{x - A} {\overline F (x - y,x
- y + 1]} dF(y) = o(F(x,x + 1]),
\end{eqnarray}
where in the second step,  Lemma \ref{jbcfby} is applied since $F\in
\mathscr{S}_{\Delta(t)}$. Substituting (\ref{erjjhfjcf1}) and
(\ref{erjjhfjcf2}) into (\ref{erjjhfj}), we obtain
\begin{eqnarray}\label{erjjhgx1}
\int_A^{x - A} {\{\overline H (x - y) - \overline H (x)\}} dF(y)& =&
K\int_A^{x - A} {\{\overline F (x - y) - \overline F (x)\}}
dF(y)\nonumber\\
&& + o(F(x,x + 1]).
\end{eqnarray}
For the same reason, in view of (\ref{GFjbdj}) and
(\ref{Gjbczsx222}), we obtain
\begin{eqnarray}\label{erjjhgx2}
\int_A^{x - A} {\{\overline H (x - y) - \overline H (x)\}} dH(y)& =&
K\int_A^{x - A} {\{\overline F (x - y) - \overline F (x)\}}
dH(y)\nonumber\\
&& + o(F(x,x + 1]).
\end{eqnarray}
By integrating by parts, we have
\begin{eqnarray}
&&\int_A^{x - A} {\{\overline F(x - y) - \overline F(x)\}} dH(y) \nonumber\\
&=& \int_A^{x - A} {\{\overline H (x - y) - \overline H (x)\}} dF(y)
+ \{\overline F (x - A) - \overline F (x)\}\{\overline H (A) -
\overline
H (x - A)\}\nonumber\\
&& - \{\overline H (x - A) - \overline H (x)\}\{\overline F (A) -
\overline F (x - A)\},
\end{eqnarray}
hence by (\ref{GFjbdj}) and Lemma \ref{gtdlyl}, we obtain
\begin{eqnarray}\label{erjjhgx3}
&&\int_A^{x - A} {\{\overline F(x - y) - \overline F(x)\}} dH(y) \nonumber\\
&=& \int_A^{x - A} {\{\overline H (x - y) - \overline H (x)\}} dF(y)
+ o(F(x,x + 1]).
\end{eqnarray}
Then from (\ref{erjjhgx1}), (\ref{erjjhgx2}), (\ref{erjjhgx3}) and
(\ref{pfgx}) it follows that
\begin{eqnarray}
&&\int_A^{x - A} {\{\overline H (x - y) - \overline H (x)\}} dH(y) -
\overline H ^2 (x) \nonumber\\
&=& K^2 \left\{ {\int_A^{x - A} {\{\overline F (x - y) - \overline F
(x)\}} dF(y) - \overline F ^2 (x)} \right\} + o(F(x,x + 1]).
\end{eqnarray}
Thus by Lemma \ref{erjdjyl} and (\ref{GFjbdj}), we conclude that
$H\in \mathscr{S}_2$. \\
{\bf Proof of Lemma \ref{lczstz}.}  Firstly, it is easy to see
\begin{eqnarray}\label{jbjf}\int_0^x F(s,s+t]ds=\int_0^t\overline{F}(s)ds-\int_x^{x+t}\overline{F}(s)ds.
\end{eqnarray}
Let $x\to\infty$ in (\ref{jbjf}), we obtain
\begin{eqnarray}\label{wqds}\int_0^{\infty} F(s,s+t]dt=\int_0^t\overline{F}(s)ds<\infty,
\end{eqnarray}
and thus $KF(\cdot+\Delta(t))$ is a density function. From Proposition 2, the
proof of Lemma 1 of Asmussen et al. \cite{asmussen2003}, it is easy to see that
$F\in \mathscr{S}_{\Delta(t)}$ is equivalent to that
 $F\in\mathscr{L}_{\Delta(t)}$ and for every function $l(x)$ such that
 $l(x)\rightarrow\infty$ and
$l(x)<x/2$, the following relation holds:
\begin{equation}\label{lx}
\int_{l(x)}^{x - l(x)} {F(x - y + \Delta(t) )} dF(y) = o(F(x + \Delta(t)
)),\begin{array}{*{20}c}
   {} & {x \to \infty .}  \\
\end{array}
\end{equation}
Note that if (\ref{lx}) holds with $l(x)$ replaced by some $l_1(x)$ such that $l_1(x)<l(x)$, then (\ref{lx}) itself holds. Hence without loss of generality, we assume $t$ divides exactly $x -
2l(x)$ and denote $n(x)=(x - 2l(x))/t$. Assume
$F\in\mathscr{L}_{\Delta(t)}$. Then we have
\begin{eqnarray}
&&\int_{l(x)}^{x - l(x)} {F(x - y + \Delta(t) )} dF(y) \nonumber\\
&=& \sum\limits_{k = 1}^n {\int_{l(x) + (k - 1)t}^{l(x) + kt} {F(x -
y+ \Delta(t) )} dF(y)}\nonumber\\
 &\sim & \sum\limits_{k = 1}^n {F(x - l(x) - (k - 1)t + \Delta(t) )F(l(x) + (k - 1)t + \Delta(t))}\nonumber\\
&\sim & \frac{1}{t}\sum\limits_{k = 1}^n {\int_{l(x) + (k - 1)t}^{l(x) + kt}
{F(x - y + \Delta(t) )} F(y + \Delta(t)
 )dy}\nonumber\\
 &=& \frac{1}{t}\int_{l(x)}^{x - l(x)} {F(x - y + \Delta(t) )} F(y + \Delta(t) )dy,\quad x \to \infty.
\end{eqnarray}
Thus, by Proposition 6 of Asmussen et al. \cite{asmussen2003}, we prove the desired result. \\
 {\bf Proof of Proposition
\ref{condition1}.} Note that
\begin{eqnarray}\label{hkz}&&\int_{-1}^{y}h(x-t-1)dt\nonumber\\
&=&
\int_{y}^{y+1}\overline{F}(x-t)dt-\int_{-1}^{0}\overline{F}(x-t)dt\nonumber\\
&\ge & \overline{F}(x-y)-\overline{F}(x),
\end{eqnarray}
hence,
\begin{eqnarray}\label{er3}&&\int_A^{x/2} {\{\overline F (x - y) - \overline F (x)\}}
dF(y)\nonumber\\
&\le & \int_A^{x/2} \int_{-1}^{y}h(x-t-1)dt
dF(y)\nonumber\\
&\le &
\overline{F}(A)\int_{-1}^{A}h(x-t-1)dt+\int_{A}^{x/2}h(x-t-1)\overline{F}(t)dt\nonumber\\
&:=& V_1+V_2,
\end{eqnarray}
where in the second step, Fubini's theorem is applied to interchange
the order of integration. It is easy to see
\begin{eqnarray}\label{V1c}V_1\sim A\overline{F}(A)h(x)=o(h(x)).
\end{eqnarray}
Since
\begin{eqnarray}\int_{0}^{x}\overline{F}(x-t)\overline{F}(t)dt=2\int_{0}^{x/2}\overline{F}(x-t)\overline{F}(t)dt,
\end{eqnarray}
it is easy to see that $F\in \mathscr{S}^*$ implies
\begin{eqnarray}\int_{A}^{x/2}\overline{F}(x-t)\overline{F}(t)dt=o(\overline{F}(x)).
\end{eqnarray}
By Theorem 2.0.8 in Bingham et al. \cite{bingham1987}, (\ref{qd})
holds locally uniformly in $(0,\infty)$. Hence
\begin{eqnarray}\label{V2c}V_2&\sim &
\int_{A}^{x/2}h(x-t)\overline{F}(t)dt\nonumber\\&=&\int_{A}^{x/2}q(x-t)\overline{F}(x-t)\overline{F}(t)dt\nonumber\\
&\lesssim &
q(x)\int_{A}^{x/2}\overline{F}(x-t)\overline{F}(t)dt=o(h(x)).
\end{eqnarray}
Combining (\ref{V1c}) and (\ref{V2c}) gives
\begin{eqnarray}\label{er1}\int_A^{x/2} {\{\overline F (x - y) - \overline F (x)\}}
dF(y)=o(h(x)).
\end{eqnarray}
Hence by integrating by parts and using $\overline F ^2 (x/2) =
o(h(x))$, we have
\begin{eqnarray}\int_{x/2}^{x-A} {\{\overline F (x - y) - \overline F
(x)\}}dF(y)=\int_A^{x/2} {\{\overline F (x - y) - \overline F
(x)\}}dF(y)+o(h(x)).
\end{eqnarray}
From this and (\ref{er1}), it follows that
\begin{eqnarray}\label{er2}\int_A^{x-A} {\{\overline F (x - y) - \overline F
(x)\}}dF(y)=o(h(x)).
\end{eqnarray}
Note that $\overline F ^2 (x)\le \overline F ^2 (x/2) = o(h(x))$,
hence by Lemma \ref{erjdjyl}, we prove (\ref{ercjjgx}). From this
and Lemma \ref{T1lxtl}, it follows that $F\in
\mathscr{S}_{\Delta(t)}$ for all $t>0$ and hence, the proof is
completed.
\\
{\bf Proof of Corollary \ref{condition2}.} By (\ref{hkz}) and
(\ref{hd}), we have
\begin{eqnarray}&&\int_A^{x/2} {\{\overline F (x - y) - \overline F (x)\}}
dF(y)\nonumber\\
&\le & \int_A^{x/2} \int_{-1}^{y}h(x-t-1)dt
dF(y)\nonumber\\
&\lesssim & h(x) \int_A^{x/2} y dF(y)=o(h(x)),
\end{eqnarray}
i.e. the relation (\ref{er1}) holds. The remaining proof is similar
to that of Proposition \ref{condition1} and we omit it.\\
{\bf Acknowledgements}\\
\mbox{}\quad This research is supported by NNSF (grant No. 10926043) in China. The author is also thankful to the anonymous referee for his/her helpful comments which results in the improvement of this paper.

\end{document}